\documentclass{amsproc}
\usepackage{graphicx}      
\usepackage[all]{xy}
\usepackage{amssymb}

\issueinfo{00}
  {}
  {}
  {1995}

\copyrightinfo{1995}
  {American Mathematical Society}

\newtheorem{thm}{Theorem}[section]
\newtheorem{lm}[thm]{Lemma}
\newtheorem{cor}[thm]{Corollary}

\theoremstyle{definition}
\newtheorem{defn}[thm]{Definition}
\newtheorem{ex}[thm]{Example}

\theoremstyle{remark}
\newtheorem{re}[thm]{Remark}

\numberwithin{equation}{section}

\begin{document}

\title{A class of gradings of simple Lie algebras}

\author{Karin Baur}
\address{Department of Mathematics, 
University of Leicester, University Road, 
Leicester LE1 7HR, United Kingdom}
\email{k.baur@mcs.le.ac.uk}
\thanks{The first author was supported by the 
Freie Akademische Stiftung and 
by DARPA contract \#AFRL F49620-02-C-0010}

\author{Nolan Wallach}
\address{Department of Mathematics, University 
of California at San Diego, La Jolla, CA 92122}
\email{nwallach@ucsd.edu}

\subjclass{Primary ; Secondary }

\date{February 20, 2006}

\keywords{Parabolic subgroups, moment map}

\begin{abstract}
In this paper we give a classification  
of parabolic subalgebras
of simple Lie algebras over $\mathbb{C}$ that satisfy 
two properties. The
first property is Lynch's sufficient condition for the 
vanishing of certain
Lie algebra cohomology spaces for generalized Whittaker 
modules associated
with the parabolic subalgebra and the second is that the 
moment map of the
cotangent bundle of the corresponding generalized flag variety 
be birational
onto its image. We will call this condition the moment map 
condition. 
\end{abstract}

\maketitle
%
%

%
\section{Introduction}\label{sec:intro} 
%

The purpose of this paper is to give a classification of parabolic 
subalgebras
of simple Lie algebras over $\mathbb{C}$ that satisfy 
two properties. The
first property is Lynch's sufficient condition for the 
vanishing of certain
Lie algebra cohomology spaces for generalized Whittaker 
modules associated
with the parabolic subalgebra and the second is that the 
moment map of the
cotangent bundle of the corresponding generalized flag variety 
be birational
onto its image. We will call this condition the moment map 
condition. 
Associated to each parabolic subalgebra of a simple Lie algebra $\mathfrak{g}$ 
is a $\Bbb{Z}$-grading and to each $\Bbb{Z}$-grading of $\mathfrak{g}$ 
corresponds a parabolic subalgebra. The first condition is that the parabolic 
subalgebra has a Richardson element in the first graded part (where the grading 
is the one associated to the parabolic subalgebra).

If $G$
is a semi-simple Lie group over $\mathbb{R}$ and if $P$ is a 
parabolic
subgroup of $G$ such that the intersection of the 
complexification of its Lie
algebra intersected with each simple factor of the 
complexification of
$Lie(G)$ satifies the two conditions then one can prove holomorphic
continuation of Jacquet integrals and a variant of a multiplicity 
one theorem
for degenerate principal series associated with $P$ 
(\cite{w}, \cite{y}). The full
classification corresponding to the first condition was the 
subject of our
joint paper~\cite{bw} where we classified the 
so-called \textquotedblleft
nice\textquotedblright\ parabolic subalgebras of simple 
Lie algebras over
$\mathbb{C}$. Thus the point of this paper is to list the elements 
of the list
in~\cite{bw} that satisfy the second condition. The second condition 
is just the
assertion that the stabilizer of a Richardson element in the 
nilradical of the
parabolic subgroup corresponding to the parabolic subalgebra of 
the adjoint
group  is the same as the stabilizer in the adjoint group. In the 
case of type
$A_{n}$ this condition is automatic (cf.~\cite{h}) thus nice 
implies both
conditions. For the other classical groups the moment map condition 
is not
automatically satisfied. One especially striking aspect of this 
classification
is that the exceptional groups are \textquotedblleft better
behaved\textquotedblright\ than the classical groups. Indeed, 
except for the
case of $E_{7}$ every nice parabolic subalgebra of an exceptional 
simple Lie
algebra also satisfies the moment map condition. 
In the case of $E_{7}$ there
are $3$ (out of the $29$) that do not.

One can argue that the results of this paper are for all practical 
purposes in
the literature. Indeed, most of our work involves the classical 
groups and
amounts to explaining how certain results of Hesselink~\cite{h} 
(which he basically
attributes to~\cite{sp}) apply. However, if one is not an expert 
in the subject
then the determination of whether a parabolic subgroup satisfies the
conditions of the holomorphic continuation and multiplicity 
one theorem would
involve a serious effort in a field only peripherally related 
to the desired
application. We therefore felt that there was value in placing 
the necessary
results in a short paper that could be understood with only the 
knowledge of
the essentials of algebraic group theory and to give an easily 
applied listing
of the pertinent parabolic subalgebras.

The paper is organized as follows: 
The full classification is described in the next section. In section 3 
we discuss the moment map condition and consider the stabilizer 
subgroup in $P$ of a Richardson element. Section 4 treats the proof 
of the theorem in the classical case. Following that, the proofs for 
the exceptional cases are found in section 5. 
In the appendix we list all parabolic subalgebras of the exceptional 
Lie algebras where there is a Richardson element in the first graded 
part and such that the moment map condition is satisfied. 

The authors would like to thank Hanspeter Kraft for helpful discussions 
concerning the moment map condition and the normality.

%
\section{Statements of the results}\label{sec:one} 
%
If not specified otherwise, $\mathfrak{g}$ will denote a simple Lie 
algebra over the complex numbers. 
Fix a Borel subalgebra $\mathfrak{b}$ in $\mathfrak{g}$, let 
$\mathfrak{h}\subset\mathfrak{b}$ be a Cartan subalgebra of 
$\mathfrak{g}$. We 
will denote the set of simple roots relative to this choice by 
$\Delta=\{\alpha_1,\dots,\alpha_n\}$. We always use the 
Bourbaki-numbering of simple roots. 

Let $\mathfrak{p}\subset\mathfrak{g}$ be a parabolic subalgebra, 
$\mathfrak{p}=\mathfrak{m}\oplus\mathfrak{u}$ (where $\mathfrak{m}$ 
is a Levi 
factor and $\mathfrak{u}$ the corresponding nilpotent radical of 
$\mathfrak{p}$). After conjugation we can assume that $\mathfrak{p}$ 
contains the chosen Borel subalgebra and 
$\mathfrak{m}\supset\mathfrak{h}$. If $\mathfrak{b}$ has been fixed 
then we 
will say that $\mathfrak{p}$ is standard if  
$\mathfrak{p}\supset\mathfrak{b}$ 
from now on. In particular, if $\mathfrak{p}$ is standard then it is 
given by a subset of $\Delta$, namely the simple roots such that 
both roots spaces $\mathfrak{g}_{\pm\alpha}$ belong to the Levi 
factor 
of $\mathfrak{p}$. 
 
Thus such a parabolic subalgebra is described by a $n$-tuple, 
$(u_1,...,u_n)$ in $\{0,1\}^n$: ones correspond to simple roots 
with root spaces not in $\mathfrak{m}$. Equivalently, a parabolic 
subalgebra is given by a coloring of the Dynkin diagram of the Lie 
algebra: a black (colored) node corresponds to a simple root whose 
root space belongs to $\mathfrak{m}$. Here, one has to be very 
careful 
since there exist different notations. Our choice was motivated by 
the coloring for Satake diagrams.  Let $(u_1,...,u_n)$ define the 
parabolic subalgebra $\mathfrak{p}$ and and $H \in \mathfrak{h}$ be 
defined by $\alpha_i(H) = u_i$. If we set $\mathfrak{g}_i=\{x\in 
\mathfrak{g}|[H,x]=ix\}$ then 
$\mathfrak{p}=\sum_{i\geq 0}\mathfrak{g}_i$. 
 
%
%
\subsection{Results in the classical cases} 
%
%
As is usual, we will refer to the simple Lie algebras of type 
${\rm A}_n,{\rm B}_n,{\rm C}_n,{\rm D}_n$ as the classical Lie 
algebras and the remaining five simple Lie algebras will be called 
exceptional.  We realize the classical Lie algebras as subalgebras 
of $\mathfrak{gl}_N$ for $N=n+1, 2n+1,2n,2n$ respectively. With 
${\rm A}_n$ the trace zero matrices, ${\rm B}_n,{\rm D}_n$ the 
orthogonal Lie algebra of the symmetric form with matrix with all 
entries $0$ except for those on the skew diagonal which are $1$ 
and ${\rm C}_n$ the symplectic Lie algebra for the symplectic 
form with matrix whose only nonzero entries are skew diagonal and 
the first $n$ are $1$ and the last $n$ are $-1$. With this 
realization we take as our choice of Borel subalgebra the 
intersection of the corresponding Lie algebra with the upper 
triangular matrices in $\mathfrak{gl}_N$. We will call a parabolic 
subalgebra that contains this Borel subalgebra standard. If 
$\mathfrak{p}$ is a standard parabolic subalgebra then we refer 
to the 
Levi factor that contains the diagonal Cartan subalgebra, 
$\mathfrak{h}$, by $\mathfrak{m}$ and call it the standard Levi 
factor. 
One has to be careful in the case of ${\rm D}_n$ since the two 
parabolic subalgebras with exactly one of the last two simple 
roots $\alpha_n$, $\alpha_{n-1}$ as roots of the Levi factor 
are conjugate under an outer automorphism. 
To avoid ambiguity for $n\ge 5$, we will always assume in 
this case that $\alpha_{n-1}$ is the root of the Levi factor.

Thus for all classical Lie algebras the standard Levi factor is then 
in diagonal block form given by a sequence of square matrices on 
the diagonal. We denote it by $\underline{d}=(d_1,\dots,d_m)$. 
For the orthogonal and symplectic Lie algebras, 
these sequences are palindromic.  
In that case we write $\underline{d}=(d_1,\dots,d_r,d_r,\dots,d_1)$ 
or $\underline{d}=(d_1,\dots,d_r,d_{r+1},d_r,\dots,d_1)$. 
If $\mathfrak{p}$ is a parabolic 
subalgebra for one of these Lie algebras then if $\mathfrak{m}$ 
is the 
standard Levi factor of the parabolic subalgebra to which it is 
conjugate then we will say that $\mathfrak{m}$ is the standard Levi 
factor. 
 
We describe now the parabolic subgroups in the classical 
cases that have the desired properties. 

\begin{thm}\label{thm:classical}
Let $P\subset G$ be a parabolic 
subgroup of a classical group. 
Let $\mathfrak{m}$ be the standard 
Levi factor of the Lie algebra of $P$, 
let $\underline{d}$ be the vector describing the 
lengths of the blocks in $\mathfrak{m}$. 
Then there is a Richardson element $x\in \mathfrak{g}_1$ 
with $G_x=P_x$ if and 
only if $P\subset G$ is one of the following: 

\begin{itemize}
\item[(i)]
$G={\rm SL}_{n+1}$ and the entries of $\underline{d}$ satisfy 
$d_1\le\dots\le d_s\ge\dots\ge d_r$ 
({\bf unimodality} of $\underline{d}$).
\item[(ii)]
$G={\rm Sp}_{2n}$, $\underline{d}$ is unimodal and if 
$\mathfrak{m}$ has an odd number of blocks, 
all $d_i$ are even.
\item[(iii)]
$G={\rm SO}_N$, $\underline{d}$ is unimodal and if 
$\mathfrak{m}$ has an even number of 
blocks then 
there is at most one odd $d_i$ and 
in this case, $i<r$ and $d_i\le d_r-3$.
\end{itemize}
\end{thm}

As a consequence, any such $P$ is given by 
an $\mathfrak{sl}_2$-triple 
(see the next section for the definition) 
if $G={\rm Sp}_{2n}$, 
$G={\rm SO}_{2n+1}$ or in the case of $G={\rm SO}_{2n}$ with 
$\mathfrak{m}$ having an odd number of blocks. 
If $G={\rm SO}_{2n}$ with an even number of blocks, 
the only such $P$ that are not given by an 
an $\mathfrak{sl}_2$-triple are the ones with 
one odd $d_i$ (and $d_i\le d_r-3$). 

It is often useful to know when the 
orbit of Richardson elements has normal closure. 
We keep the notation of Theorem~\ref{thm:classical} above. 
By~\cite{kp0}, any orbit closure is normal for $G={\rm SL}_n$. 
Using~\cite{kp} (for all cases except for part of 
the very even orbits in $SO_{4n}$) and~\cite{so2} (for 
the remaining very even orbits), one can check whether 
the orbit of Richardson elements has normal closure. 

\begin{re}
Let $P=P(d)\subset G$ with $G={\rm Sp_{2n}}$ or ${\rm SO_N}$ 
be one of the parabolic subgroups of 
Theorem~\ref{thm:classical}. Let $X\in\mathfrak{n}$ be a Richardson 
element for $P$. Then the orbit $GX$ 
has normal closure exactly in the following cases. 

(i) $G={\rm Sp}_{2n}$ 
$
\left\{
\begin{array}{l} 
\mathfrak{m} \text{ has 
an even number of blocks.} \\ 
\mathfrak{m} \text{ has an odd number 
of blocks and } d_1=d_2=\dots=d_r.
\end{array}
\right.$
 
(ii) $G={\rm SO}_{N}$ 
$
\left\{
\begin{array}{cl} 
1) & \mathfrak{m} \text{ has an odd number of blocks.} \\
2) & \mathfrak{m} \text{ has an even number of blocks, and} \\ 
  & \text{a) } d_1=\dots=d_r \text{ (even)},\\ 
 & \text{b) } d_1=\dots=d_s,\ d_{s+1}=d_s+2,\ 
d_{s+1}=\dots=d_r, \text{ (even)},\\
 & \text{c) there is one odd entry } d_i, d_i\le d_r-3 \text{ {\bf with $i=1$} } \\
 & \text{and } d_2=\dots=d_r 
 \text{ \bf or $i>1$},\ d_1=\dots=d_{i-1}, 
 d_i-d_{i-1}=1 \\ 
 &  \text{and } d_{i+1}=\dots=d_r.
\end{array}\right.$

\end{re}
%

%
\subsection{Results in the exceptional cases} 
%
%
In this subsection we will state the classification of 
nice parabolic subalgebras for the 
exceptional simple Lie algebras. The parabolic subalgebra will be 
given by an $n$-tuple where $n$ is the rank and the entries are 
$\alpha_i(H)$ where $H$ is the element that gives the grade 
corresponding to the parabolic subalgebra and the $\alpha_i$ are 
the simple roots in the Bourbaki order. 
 
We give an explicit list of all the parabolic subgroups 
with a Richardson element in $\mathfrak{g}_1$ and 
with $P_x=G_x$ for Richardson elements in the appendix. 

We first recall that 
for $G={\rm G}_2$ or ${\rm F}_4$ all parabolic subalgebras 
with an Richardson element in 
$\mathfrak{g}_1$ are given by an $\mathfrak{sl}_2$-triple. 
(cf.~\cite{bw}). Since all such $P$ satisfy $P_x=G_x$ 
as we will see later, the list of parabolic subgroups 
with a Richardson element in $\mathfrak{g}_1$ and 
with $P_x=G_x$ is just the list of parabolic subgroups 
given by an $\mathfrak{sl}_2$-triple. 

If $G$ is of type ${\rm E}_n$ than the picture is the following: 

\begin{thm}
The parabolic subgroups $P$ of ${\rm E}_n$ that have 
a Richardson element $x$ in $\mathfrak{g}_1$ 
and $P_x=G_x$ are the parabolic subgroups 
that have a Richardson element in $\mathfrak{g}_1$ 
except for the following three subgroups of ${\rm E}_7$: 
\begin{eqnarray*}
& (1,1,0,0,0,0,1) \\ 
& (0,0,1,0,0,0,1) \\
& (0,0,0,0,1,0,1)
\end{eqnarray*}
\end{thm}
%
%
%
\section{Birationality of the moment map}
%
Let $G$ be a connected semisimple linear algebraic 
group over the 
complex numbers, $B$ a Borel subgroup containing 
the maximal torus $T$. We denote the Lie algebras 
with corresponding gothic letters.
Let $P\supset B$ be a parabolic 
subgroup of $G$, $P=M\cdot N$ with $M$ a Levi factor 
and $N$ the corresponding unipotent radical.

The dual of the cotangent bundle of $G/P$, $T^*(G/P)$ 
is a $G$-manifold with a natural symplectic form. In 
particular, there exists a moment map 
\[
\mu_P:\ T^*(G/P) \longrightarrow \mathfrak{g}^*
\]
After identifying $T^*(G/P)$ with the homogeneous 
$G$-vector bundle $G\times^P(\mathfrak{g/p})^*$ 
and dualizing (use 
$G\times^P(\mathfrak{g/p})^*\cong G\times^P\mathfrak{n}$) 
we obtain the map 

\begin{eqnarray*}
\psi_P:\ G\times^P\mathfrak{n} & \longrightarrow
 & \mathfrak{g} \\ 
 (g,X) & \mapsto & Ad(g)X 
\end{eqnarray*}
where $\mathfrak{n}=\rm{Lie}\,N$ is the nilradical of $P$. 
Its image is $\overline{Gx}$ for $x$ a Richardson element 
in $\mathfrak{n}$. 
In case $P=B$, a Richardson element is a 
regular nilpotent element and so the image is equal to the 
nullcone $G\mathfrak{n}=\mathcal N$ which is normal. 

Of special interest is the case where $\psi_P$ is birational 
onto its image. In that case, a result of the second named 
author shows that for the space of generalized Whittaker 
vectors (certain linear maps on a representation induced 
from a Fr\'echet representation space of $P_{\Bbb{R}}$) 
a multiplicity one theorem holds if $P$ has a Richardson 
element in $\mathfrak{g}_1$ (cf.~\cite{w}). 

From the discussion above we find that the 
following holds:
\begin{lm}
Let $P$ be a parabolic subgroup of $G$. Then 
the induced map $\psi_P$ is birational onto its 
image if and only if for any Richardson element 
$x$ for $P$ the stabilizer $G_x$ is contained in $P$. 
\end{lm}

To describe the main class of parabolic subgroups 
for which there is a Richardson element in $\mathfrak{g}_1$ and 
such that the stabilizers of Richardson elements in 
$G$ resp. in $P$ agree, we introduce the following notion: 

\begin{defn}
Let $\mathfrak{p}=\oplus_{j\ge 0}\mathfrak{g}_j$ 
be a parabolic subalgebra where the grading of 
$\mathfrak{g}$ is given by $H\in\mathfrak{h}$. 
We say that $P$ (or $\mathfrak{p}$) is 
{\itshape given by an $\mathfrak{sl}_2$-triple} 
if there exists a nonzero $x\in\mathfrak{g}_1$ 
such that the Jacobson-Morozov triple through 
$x$ is $\{x,2H,y\}$ (for some $y\in\mathfrak{g}$). 
\end{defn}
By construction, such an $x$ is a 
Richardson element 
for $P$ that belongs to $\mathfrak{g}_1$. 
Furthermore, we recall the following standard result 
(cf.~\cite{pre}, a simple proof can be found in~\cite{w}). 

\begin{thm}\label{th:sl2}
Let $P\subset G$ be a parabolic subgroup. 
If $P$ is given by an $\mathfrak{sl}_2$-triple 
then $P_x=G_x$ for Richardson elements. 
\end{thm}

We will see that if $P$ has a Richardson element $x$ in 
$\mathfrak{g}_1$ and $P_x=G_x$ then $P$ is given by 
an $\mathfrak{sl}_2$-triple for types ${\rm B}_n,{\rm C}_n$, 
${\rm G}_2$ and ${\rm F}_4$ and for 
${\rm D}_n$ if the number of 
blocks in $\mathfrak{m}$ is odd.

The following observations proves to be very useful 
for the exceptional cases.
\begin{lm}\label{lm:stabil-in-levi}
If $P$ is a parabolic subgroup that has a Richardson 
$x$ element in $\mathfrak{g}_1$ then the number of 
components of $P_x$ is the same as the number of 
components of $M_x$. 
\end{lm}
\begin{proof}
Let $p=n\cdot m\in P_x$ , $n\in N$ and $m\in M$. 
Set $y=Ad(m)x$. This lies in $\mathfrak{g}_1$ since 
$Ad(m)$ preserves 
every component $\mathfrak{g}_j$. 
Now $n$ is unipotent, so 
$x=Ad(p)x=Ad(n)(Ad(m)x)=Ad(n)y$ 
$=y\mod \sum_{j>1}\mathfrak{g}_j$. 
So $y=x$, i.e. $m\in M_x$ and then 
$n\in N_x$. 
This gives $P_x=M_x N_x$, hence 
$P_x/P_x^o=M_x/M_x^o$ ($N_x$ is connected). 

%
\end{proof}

\begin{lm}\label{lm:levi-factor}
Let $P_i\subset G$, $i=1,2$,  be parabolic subgroups 
with Richardson elements $x_i$ 
in $\mathfrak{g}_1$. Assume that the corresponding 
Levi factors are conjugate. 
Then $x_1$ and $x_2$ are conjugate.
\end{lm}
\begin{proof}
This is Corollary 5.18 in~\cite{bj}. 
\end{proof}

\begin{cor}\label{cor:conj-levi}
In the situation above, the stabilizer subgroups 
$(M_1)_{x_1}$ and 
$(M_2)_{x_2}$ have the same number of components. 
\end{cor}

%
\section{Proof of the theorem, classical case}
%

We prove the statements of Theorem~\ref{thm:classical} 
by a case by case check. 

Recall that if $G$ is classical case, $P\subset G$ 
is given by the sequence of block lengths in the 
standard Levi factor, $\underline{d}=(d_1,\dots,d_m)$ 
for ${\rm A}_n$ and 
$\underline{d}_{even}=(d_1,\dots,d_r,d_r,\dots,d_1)$ 
or 
$\underline{d}_{odd}=(d_1,\dots,d_r,d_{r+1},d_r,\dots,d_1)$ 
(with $d_{r+1}>0$) for ${\rm B}_n$, ${\rm C}_n$, ${\rm D}_n$. 
If $G$ is of type ${\rm A}_n$, then we always have 
$P_x=G_x$ ($x$ nilpotent), since 
$G_x=\mathfrak{gl}(n+1,\Bbb{C})_x\cap G$ or 
as Hesselink observed in~\cite[Section 3.1]{h}: 

\begin{lm}\label{lm:A-connected}
Let $x\in\mathfrak{gl}_n$ be a nilpotent element. 
Then the centralizer $G_x$ of $x$ in $G$ is connected. 
\end{lm}

In our earlier paper~\cite{bw} we have described 
the parabolic subgroups of the classical groups 
that are given by an $\mathfrak{sl}_2$-triple. 
In type ${\rm A}_n$, $P$ is given by a $\mathfrak{sl}_2$-triple 
if and only if $\underline{d}$ is unimodal and palindromic. 
On the other hand, any parabolic subgroup of ${\rm A}_n$ 
with unimodal sequence of block lengths in the 
standard Levi factor has a Richardson element 
in $\mathfrak{g}_1$. So there are many parabolic 
subgroups that are not given by an $\mathfrak{sl}_2$-triple 
but still have a Richardson element in $\mathfrak{g}_1$ 
and birational $\psi_P$.

From a more general statement of Hesselink for 
the number of components of $G_x/P_x$ for types ${\rm B}_n$, 
${\rm C}_n$ and ${\rm D}_n$ 
we obtain the characterization of parabolic subgroups of 
these with $G_x=P_x$ 
($x$ a Richardson element).

To formulate the statement we introduce some notation: 
As is customary, set $\epsilon=1$ for $G={\rm Sp}_{2n}$ 
and $\epsilon=0$ for $G={\rm SO}_N$.  
If $x$ is a Richardson element and $\lambda$ its partition, 
ordered as 
$\lambda=\lambda_1\ge \lambda_2\ge\dots$, we define 
$N_{odd}(\lambda)$ to be the number of odd parts in 
$\lambda$ and 
$B(\lambda):=\{j\in\Bbb{N}\mid \lambda_j>\lambda_{j+1}, 
\lambda_j\not\equiv\epsilon\mod 2\}$.

\begin{thm}\label{thm:hesselink}
Let $G$ be of type ${\rm B}_n$, ${\rm C}_n$ or ${\rm D}_n$, 
$P\subset G$. Let $\lambda$ be the partition of the 
Richardson element for $P$ corresponding to its Jordan form 
as an element of $M_n(\Bbb{C})$. 
Then the induced map $\psi_P$ is birational 
onto its image exactly in the following cases:
\[
\begin{array}{lllll}
(i) & N_{odd}(\lambda)=d_{r+1} & \text{if } \underline{d}
=\underline{d}_{odd} \\ 
(ii) &  N_{odd}(\lambda)=0 & \text{if }\underline{d}
=\underline{d}_{even}, 
 &  G={\rm Sp}_{2n}; \\
   & & \text{if }d=\underline{d}_{even}, &  G={\rm SO}_{2n},
 & B(\lambda)=\emptyset;\\
(iii) & N_{odd}(\lambda)=2 &\text{if } \underline{d}
=\underline{d}_{even},
 & G={\rm SO}_{2n}, & 
B(\lambda)\neq\emptyset. 
\end{array}
\]
\end{thm}
\begin{proof}
This is a special case of Theorem 7.1 in~\cite{h}. 
\end{proof}

As an immediate consequence, we obtain:
\begin{cor}\label{cor:sp-even}
Let $P$ be a parabolic subgroup of ${\rm Sp}_{2n}$, 
given by $\underline{d}=\underline{d}_{even}$. 
Then $P$ has a Richardson 
element $x\in\mathfrak{g}_1$ and $P_x=G_x$ if and 
only if $P$ has a Richardson element in $\mathfrak{g}_1$ 
if and only if $d_1\le\dots\le d_r$. 
\end{cor}
\begin{proof}
For any $P\subset {\rm Sp}_{2n}$ given by $\underline{d}_{even}$ 
the partition of a Richardson element has only even 
parts. 
Explicitly, $\lambda=(2r)^{d_1}, (2r-2)^{d_2-d_1},\dots, 
4^{d_{r-1}-d_{r-2}}, 2^{d_r-d_{r-1}}$ 
(\cite[Lemma 4.7]{b05}). So by part (i) of 
Theorem~\ref{thm:hesselink}, $G_x=P_x$ for Richardson 
elements. 
The statement then follows since $P$ has a Richardson 
element in $\mathfrak{g}_1$ if and only if $\underline{d}_{even}$ 
is unimodal (\cite[Theorem 1.3]{bw}).
\end{proof}

\begin{cor}\label{cor:sp-odd}
Let $P\subset G={\rm Sp}_{2n}$ be given by $\underline{d}_{odd}$. 
Then $P$ has a Richardson element $x$ in $\mathfrak{g}_1$ 
and $P_x=G_x$ if and only if all $d_i$ are even.
\end{cor}

\begin{proof}
From~\cite[Theorem 1.3]{bw} we know that $P$ has a 
Richardson element 
in $\mathfrak{g}_1$ if and only if $\underline{d}_{odd}$ is 
unimodal and if any odd $d_i$ appears only once 
among $d_1,\dots,d_r$. 
From~\cite[Lemma 4.7]{b05} the dual of the partition of a 
Richardson element is 
$d_{r+1}\cup\{\bigcup_{d_i\notin D_o}d_i,d_i\}$
$\cup\{\bigcup_{d_i\in D_o}d_i-1,d_i+1\}$ 
where $D_o:=\{d_i\mid d_i\equiv 1\}$ is the set 
of odd entries of $\underline{d}_{odd}$. 
In any case, the partition of a Richardson element 
has $d_{r+1}$ parts. If all $d_i$ are even, 
the partition of a Richardson element is 
\[
\lambda: (2r+1)^{d_1},(2r-1)^{d_2-d_1},\dots,3^{d_r-d_{r-1}},
1^{d_{r+1}-d_r}
\]
In particular, its parts are all odd and by (i) of 
Theorem~\ref{thm:hesselink}, we get $G_x=P_x$ for 
Richardson elements. 

If we assume that $D_o$ is not empty (i.e. that 
there are odd $d_i$ then the partition $\lambda$ 
of a Richardson 
element still has $d_{r+1}$ parts. As one can check, 
$\lambda$ contains $2$ even parts for every odd 
$d_i$. So by (ii) of Theorem~\ref{thm:hesselink}, 
$\psi_P$ is not birational (it is in fact a 
$k$-fold covering map where $k=2^{d_{r+1}/2-|D_o|}$, 
using the full version of Theorem 7.1 of~\cite{h}). 
\end{proof}

\begin{cor}\label{cor:so-even}
Let $P$ be a parabolic subgroup of ${\rm SO}_{2n}$ 
that is given by $\underline{d}=\underline{d}_{even}$. 
Then $P$ 
has a Richardson element $x$ in $\mathfrak{g}_1$ 
and $G_x=P_x$ if and only there is at most 
one odd $d_i$ and in that case, $d_i\le d_r-3$. 
\end{cor}

\begin{proof}
We recall Theorem 1.4 of~\cite{bw}: $P$ has a 
Richardson element in $\mathfrak{g}_1$ 
exactly in the following cases: 
$d_1\le\dots\le d_r$ or 
$d_1\le\dots\le d_{t-1}<d_t$ and $d_{t+1}=\dots=d_r$ are 
equal to $d_t-1$ (for $t=1$ 
the condition $d_{t-1}<d_t$ means $d_t>1$). Furthermore, 
in both cases, odd $d_i$ appear only once among 
$d_1,\dots,d_r$. 

In the second case we may reorder $\underline{d}$ 
to have $d_1\le\dots\le d_r$. The corresponding 
parabolic subalgebra still has a Richardson element 
in $\mathfrak{g}_1$. Since the two Levi factors 
are conjugate, we can apply Corollary~\ref{cor:conj-levi} 
to see 
that the stabilizers of the corresponding Richardson 
elements have the same number of components. 
Furthermore, the corresponding Richardson elements are 
conjugate. So using Lemma~\ref{lm:stabil-in-levi}, 
we see that $G_x=P_x$ can be checked for the 
reordered $\underline{d}$. 

So from now on we assume $d_1\le\dots\le d_r$. 
The dual of the 
partition of a Richardson element is 
$\{\bigcup_{d_i\notin D_o} d_i,d_i\}\cup
\{\bigcup_{d_i\in D_o} d_i-1,d_i+1\}$, where 
$D_o$ is the set of odd $d_i$ among $d_1,\dots,d_r$
(\cite[Lemma 4.7]{b05}). 
The corresponding partition 
has no odd parts if all $d_i$ are even, in that case, 
it is 
\[
2r^{d_1}, 2(r-1)^{d_2-d_1},\dots, 2^{d_r-d_{r-1}}.
\]
In particular, we then have 
$B(\lambda)=\emptyset$,  and 
by (ii) of Theorem~\ref{thm:hesselink}, $G_x=P_x$. 
Now let $D_0\neq\emptyset$. For every 
element in $D_o$, the partition of 
a Richardson element has two odd entries, i.e. 
$N_{odd}(\lambda)=2|D_o|$ and by (ii) and (iii), 
$|D_o|\le 1$. 
It remains to consider the case $D_o=\{d_i\}$. 
One checks that the corresponding partition is 
\[
\begin{array}{l}
2r^{d_1},\dots,2(r-(i-2))^{d_{i-1}-d_{i-2}},
2(r-(i-1))^{d_i-d_{i-1}-1},
(2(r-(i-1))-1)^2, \\ 
2(r-i)^{d_{i+1}-d_i-1}, 
2(r-(i+1))^{d_{i+2}-d_{i+1}},\dots, 4^{d_{r-1}-d_{r-2}},
2^{d_r-d_{r-1}}.
\end{array}
\]
It has two odd entries. 
If $d_i=d_r-1$ then $d_i+1=d_{i+1}=\dots=d_r$ and the smallest 
nonzero entries of the partition are $2(r-(i-1))-1$, i.e. odd. 
So $B(\lambda)$ is empty and hence by (ii) of 
Theorem~\ref{thm:hesselink}, $G_x\neq P_x$. 
Otherwise, $d_i\le d_r-3$ and the smallest 
entries are even, so $B(\lambda)\neq\emptyset$ and 
$G_x=P_x$.

The unimodality now follows with the observation that 
there is no odd $d_i$ with $|\max_j\{d_j\}-d_i|=1$. 
\end{proof}
It remains to deal with the case 
of the special orthogonal groups in the case where 
the standard 
Levi factor has an odd number of blocks. By 
Theorem~\ref{thm:hesselink}, the induced map 
$\psi_P$ is birational onto its image if and 
only if the partition of a Richardson element 
has exactly $d_{r+1}$ odd parts. 
\begin{cor}\label{cor:so-odd}
Let $P\subset G={\rm SO}_N$ be given by $\underline{d}
=\underline{d}_{odd}$. 
Then $P$ has a Richardson element $x$ in $\mathfrak{g}_1$ 
and $G_x=P_x$ if and only if 
$d_1\le\dots\le d_r\le d_{r+1}$. 
\end{cor}

\begin{proof}
We recall that $P$ has a Richardson element in $\mathfrak{g}_1$ 
in exactly the following cases: 
either $d_1\le\dots\le d_{r+1}$ or there is a $t\le r$ 
such that $d_1\le\dots\le d_{t-1}<d_t$ and 
$d_{t+1}=\dots=d_{r+1}$ is equal to $d_t+1$. 

We first consider the case where $\underline{d}_{odd}$ is 
unimodal (so $d_{r+1}$ is maximal). By Lemma 4.7 of~\cite{b05}, 
the dual of the partition of a Richardson element 
has just the entries $d_1,d_1,\dots,d_r,d_r,d_{r+1}$ 
and so the partition $\lambda$ of a Richardson 
element is
\[
(2r+1)^{d_1},(2r-1)^{d_2-d_1},\dots 3^{d_r-d_{r-1}},
1^{d_{r+1}-d_r}
\]
and has $d_{r+1}$ odd parts. 
Thus by (i) of Theorem~\ref{thm:hesselink}, the map 
$\psi_P$ is birational onto its image.

Now we consider the case where $\underline{d}_{odd}$ is not 
unimodal. 
As in the proof of Corollary~\ref{cor:so-even} we 
can reorder $\underline{d}$ since the corresponding 
parabolic subalgebra also has a Richardson element 
in $\mathfrak{g}_1$. So we may assume 
$d_1\le\dots\le d_r$ $d_{r+1}=d_r+1$. 

By the Lemma~\ref{lm:parts-d-odd} below, 
$\lambda$ has $d_{r+1}+2$ entries. If we show 
that all of these are odd, then by (i) of 
Theorem~\ref{thm:hesselink}, the induced map 
$\psi_P$ is not birational onto its image. 
In fact, it will be a two-fold cover. 
Let $\underline{d'}=\underline{d'}_{odd}$ be the dimension 
vector obtained 
by replacing $d_r$ by $d_r-1$. Then $\underline{d'}$ is unimodal. 
Let $X(d)$ resp. $X(d')$ be the corresponding 
Richardson elements. It is easy to see that the maps 
$X(d)^k$ and $X(d')^k$ have the same rank for every 
$k\in\Bbb{N}$ (for details we refer the reader 
to~\cite[Section 3]{bw}). Hence 
$\dim\ker X(d)^k=\dim\ker X(d')^k+2$ for all $k$. 
From that one can compute the partitions $\lambda$ 
resp. $\lambda'$: $\lambda$ is obtained from 
$\lambda'$ by adding $\{1,1\}$. We already 
know that the partition $\lambda'$ has $d_{r+1}$ 
odd entries. Hence $\lambda$ has $d_{r+1}+2$ 
odd entries. 
\end{proof}

\begin{lm}\label{lm:parts-d-odd}
Let $P\subset G$ be a parabolic subgroup of 
${\rm SO}_N$ or ${\rm Sp}_{2n}$, $P$ be given 
by $\underline{d}=\underline{d}_{odd}$. 
If max $d_i=d_{r+1}+1$ then 
the partition $\lambda$ of a Richardson element 
has $d_{r+1}+2$ parts.
\end{lm}

\begin{proof}
For $X\in\mathfrak{gl}_N$ is a nilpotent element 
with partition $\lambda$, the number of parts 
of $\lambda$ is equal to the dimension of the 
kernel of the map $X$. This follows immediately 
from the formula for the partition for nilpotent 
matrices, as given in~\cite[Section 3]{bw}: 
\[
\lambda: m^{a_m},(m-1)^{a_{m-1}},\dots, 1^{a_1}
\]
with 
$a_j:=2\dim\ker X^j-\dim\ker X^{j-1}-\dim\ker X^{j+1}$, 
$m$ the maximal number such that 
$\ker X^{m-1}\subsetneq \ker X^m=\Bbb{C}^N$ 
(with $\dim\ker X^{m+1}=\Bbb{C}^N$). In particular, 
$\lambda$ has $\sum_{k=1}^m a_k=\dim\ker X$ parts. 

Now if $X$ is a Richardson element for $P$ it is 
in particular a generic element of the nilradical. 
W.l.o.g. we can assume $d_1\le\dots\le d_r$ 
(by Lemma~\ref{lm:levi-factor}). 
So the rank of $X$ is 
\begin{eqnarray*}
{\rm rk}\,X & = & 2\sum_{i=1}^{r-1}\min\{d_1,d_{i+1}\} 
 + 2\min\{d_r,d_{r+1}\} \\ 
 & = & 2(d_1+\dots+d_{r-1}) + 2d_{r+1}
\end{eqnarray*}
Hence its kernel has dimension 
$2d_r-d_{r+1}$ 
$=d_{r+1}+2$.
\end{proof}

From the description of parabolic subalgebras given 
by an $\mathfrak{sl}_2$-triple that is in our earlier 
paper~\cite{bw} we can now deduce:
\begin{cor}
For types ${\rm B}_n, {\rm C}_n$ and ${\rm D}_n$, 
the parabolic subgroups with a 
Richardson 
element in $\mathfrak{g}_1$ and with $P_x=G_x$ are those 
given by an $\mathfrak{sl}_2$-triple except in the case 
of ${\rm D}_n$ with an even number of blocks. 
There, the extra family of examples are the ones with one 
odd entry $d_i$ (with $d_i\le d_r-3$). 
\end{cor}

\begin{re}
The proofs of the statements Corollary~\ref{cor:sp-even}, 
\ref{cor:sp-odd},
\ref{cor:so-even} and \ref{cor:so-odd} all directly used 
Theorem~\ref{thm:hesselink} via the partition of the 
corresponding Richardson element. The statements could 
have also been proved indirectly using 
Corollary 7.7 of~\cite{h}. Since this corollary uses 
very intricate notation and is stated without proof, 
we decided to include the proofs for our statements. 
\end{re}

%
\section{Result for exceptional groups}
%

The parabolic subgroups of ${\rm G}_2$, ${\rm F}_4$ and 
${\rm E}_6$ with have a Richardson element $x$ in 
$\mathfrak{g}_1$ all satisfy $G_x=P_x$. 
In the case of ${\rm G}_2$ and ${\rm F}_4$ 
this follows from the fact that all parabolic 
subgroups with a Richardson element in $\mathfrak{g}_1$ 
are given by an $\mathfrak{sl}_2$-triple. 
Hence the claim follows by Theorem~\ref{th:sl2}. 

For ${\rm E}_n$ we can use a list of Sommers. 
In~\cite{so}, Sommers classified the nilpotent elements 
of the exceptional Lie groups whose centralizers 
are not connected. For ${\rm E}_6$ all these 
are the nilpositive element of an $\mathfrak{sl}_2$-triple. 
In particular, they all have the property $G_x=P_x$ 
for Richardson elements. 

It remains to understand ${\rm E}_7$ and ${\rm E}_8$. 
From the classification of parabolic subgroups with 
a Richardson element in $\mathfrak{g}_1$ we know that 
there are five parabolic subgroups in ${\rm E}_7$ resp. 
one in ${\rm E}_8$ that are not given by an 
$\mathfrak{sl}_2$-triple but still have a Richardson 
element in $\mathfrak{g}_1$. We list them here 
together with the dimension of the Richardson orbit 
$\mathcal{O}_x$ 
(obtained from the formula 
$\dim\mathfrak{m}=\dim\mathfrak{g}^x$ that holds for 
Richardson elements) and its Bala-Carter label.
Recall that an entry $1$ stands for a simple root that is 
not a root of the standard Levi factor. 
\[
\begin{array}{cclll}
&  & & \dim\mathcal{O}_x & \text{Label}\\
 & & \\
{\rm E}_7\quad & a) & (1,1,0,0,1,0,1) & 118 & {\rm D}_6 \\ 
    & b) & (1,1,0,0,0,0,1) & 106 & {\rm D}_5(a_1)\\ 
    & c) & (0,1,1,0,0,1,1) & 118 & {\rm D}_6 \\ 
    & d) & (0,0,1,0,0,0,1) & 104 & {\rm A}_4+{\rm A}_1\\
    & e) & (0,0,0,0,1,0,1) & 104 & {\rm A}_4+{\rm A}_1\\
    &  & \\
{\rm E}_8\quad & f) & (0,0,1,0,0,0,1,0) & 216 & {\rm D}_6
\end{array}
\]
Note that the Levi factors of a) and c) are the 
same as well as the Levi factors of e) and d). 
Since the condition $P_x=G_x$ 
only depends on the Levi factor 
(cf. Lemma~\ref{lm:levi-factor}), it is 
enough to understand one of e) and d) resp. one 
of a) and c). 
Using Sommers list one checks that the groups $G_x$ 
are connected for a), c) and f). So in these cases, 
the parabolic subgroup has the property $G_x=P_x$. 

\begin{lm}
Let $P$ be one of the remaining parabolic subgroups 
b),d), e), with $x$ a Richardson element for $P$. 
Then $G_x\supsetneq P_x$.
\end{lm}

We discuss b) and d) to prove this claim. 

\begin{ex}
Let $P$ be the parabolic subalgebra 
(0,0,1,0,0,0,1) of ${\rm E}_7$. 
Set $X$ to be the following element of 
$\mathfrak{g}_1$: 
\[
X=x_3+x_{24567}+x_{134}+x_{234^256} + x_{12345}
\]
where $x_{I}$ denotes a non-zero element of the root 
space of $\alpha_I:=\sum_{i\in I}\alpha_i$. 
Using the program GAP for the Lie algebra 
$L:={\rm E}_7$ one computes that the dimension 
of the centralizer of $X$ in $L$ is $29$, 
the corresponding command is

$LieCentralizer(L, Subalgebra(L,[X]));$

\noindent
But this is just the dimension of the Levi 
factor, so $X$ is a Richardson element for 
this $P$. 

From Sommers list~\cite{so} we know that $G_X$ 
is disconnected. 
We now show that $P_X$ is connected. 
Lemma~\ref{lm:stabil-in-levi} implies that it is enough 
to prove 
that $M_x$ is connected. We will 
now calculate the stabilizer 
of $X$ in the Levi factor. If $M$ is the 
Levi factor of $P$ then $M_0:=(M,M)=$ 
${\rm GL}_2\times{\rm GL}_5$. The first graded 
part $\mathfrak{g}_1$ corresponds to the 
representation 
$(\Bbb{C}^2)^*\otimes\Lambda^2(\Bbb{C}^5)\oplus(\Bbb{C}^5)^*$. 
Under this correspondence, the Richardson element 
has the form
\[
V:=v_2\otimes(e_2\wedge e_4+e_3\wedge e_5)
+v_1\otimes(e_4\wedge e_5+e_1\wedge e_2)+e_5^*.
\]
(The $e_i$ form a basis of $\Bbb{C}^5$, the $v_i$ 
of $\Bbb{C}^*$ and $e_i^*$ of $(\Bbb{C}^5)^*$).  
Elements of $M_0$ stabilizing $V$ are 
pairs $(h,g)\in{\rm GL}_2\times{\rm GL}_5$ 
that satisfy in particular $ge_5^*=e_5^*$. 
Hence for $g=(g_{ij})_{ij}$ we have 
$g_{5i}=0$ ($i\le 4$) and $g_{55}=1$. 
Let 
\begin{small}
$h=\begin{bmatrix} a&b \\ c&d\end{bmatrix}$
\end{small}
in ${\rm GL}_2$. 
Writing out the equations given by 
$(h,g)V=V$ one uses $\det h\neq 0$ to 
obtain $g_{13}=g_{23}=g_{14}=g_{24}=0$. 
With these simplifications, the stabilizer 
in $M_0$ can be computed with mathematica, 
it is the subgroup 
\[
\left\{
\left(
\begin{bmatrix} a&0 \\ c&d\end{bmatrix},
\begin{bmatrix} 
\frac{d}{a^2} & 0 & 0 & 0 & 0 \\
g_{21} & \frac{a}{d} & 0 & 0 & -\frac{c}{d} \\ 
-\frac{c^2}{a^2d} & 0 & \frac{1}{d} & -\frac{c}{ad} & 0 \\
\frac{2c}{a^2} & 0 & 0 & \frac{1}{a} & 0 \\ 
0 & 0 & 0 & 0 & 1
\end{bmatrix} 
\right)\mid a,c,d,g_{21}\in\Bbb{C},\ ad\neq 0
\right\} \ \ \subset\ {\rm GL}_2\times{\rm GL}_5. 
\]
But this is just $(\Bbb{C}^*)^2\times \Bbb{C}^2$ 
which is connected. 
\end{ex}

\begin{ex}\label{ex:E7-1}
Let $P$ be the parabolic subalgebra 
(1,1,0,0,0,0,1) of ${\rm E}_7$. As in the previous 
example, one can use GAP to check that the 
following element 
\[
X=x_7 + x_{13456} + x_{234^25} + x_{23456} 
 - x_{234} -x_{245}
\] 
of $\mathfrak{g}_1$ satisfies 
$\dim\mathfrak{g}^X=\dim\mathfrak{m}=27$, 
so it is a Richardson element for $P$. 

By Sommers list, we know that $G_X$ is not 
connected. To see whether $P_x$ is connected 
we calculate the stabilizer of $X$ in the 
standard Levi factor $M$ of $P$. We have 
$M_0=(M,M)={\rm GL}_5$ and $\mathfrak{g}_1$ 
corresponds to the representation 
$\Bbb{C}^5\oplus\Lambda^2(\Bbb{C}^5)\oplus(\Bbb{C}^5)^*$, 
and the Richardson element has the form 
\[
V:=
e_1\wedge e_3+ e_5^* + e_1 + e_3\wedge e_5 + e_2\wedge e_4.
\]
So we have to determine whether the set 
$\{g\in{\rm GL}_5\mid gV=V\}$ is connected. 
First we observe that $ge_1=e_1$ and $ge_5^*=e_5^*$ 
imply $g_{11}=g_{55}=1$ whereas the other entries of the 
first column and last row are zero. 
Using mathematica one obtains 
\[
{\rm GL}_5\ \ni\ 
g = \begin{bmatrix}
1 & 0 & 0 & 0 & 0 \\
0 & g_{22} & 0 & g_{24} & g_{25} \\
0 & g_{32} & 1 & g_{34} & g_{35} \\ 
0 & g_{42} & 0 & g_{44} & g_{45} \\
0 & 0 & 0 & 0 & 1
\end{bmatrix}
\]
together with the equations 
\begin{eqnarray*}
g_{22}g_{34}-g_{32}g_{24}-g_{25} & = & 0 \\
g_{22}g_{44}-g_{42}g_{24} -1 & = & 0 \\
g_{32}g_{44}-g_{42}g_{34} + g_{45} & = & 0
\end{eqnarray*}
The second equation implies $g\in{\rm SL}_5$. 
We use mathematica to compute a Gr\"obner basis 
for these equations in the nine variables and 
obtain 
\[
\left\{
\begin{array}{l}
-g_{34}g_{42} + g_{32}g_{44} + g_{45}, \\
g_{34} - g_{25}g_{44}+ g_{24}g_{45}, \\
g_{32} - g_{25}g_{42} + g_{22}g_{45}, \\
-1 - g_{24}g_{42} + g_{22}g_{44}, \\ 
-g_{25} - g_{24}g_{32} + g_{22}g_{34}
\end{array}
\right\}
\]

We set $x_1=g_{22}$, $x_2=g_{44}$, $x_3=g_{42}$, $x_4=g_{24}$ 
and 
$y_1=g_{25}$, $y_2=g_{45}$ $y_3=g_{32}$, $y_4=g_{34}$. We also 
set 
\[
A=\begin{bmatrix} 
-x_4 & x_1 \\ 
-x_2 & x_3
\end{bmatrix}
\]
then the fourth equation says that $\det A=1$. Hence
\[
A^{-1}=\begin{bmatrix} 
x_3 & -x_1 \\
x_2 & -x_4
\end{bmatrix}
\]
The other four equations say that 
\[
\begin{bmatrix} 
0 & A \\ 
A^{-1} & 0
\end{bmatrix}
\begin{bmatrix} 
y_1 \\ y_2 \\ y_3 \\ y_4
\end{bmatrix}
=\begin{bmatrix} 
y_1 \\ y_2 \\ y_3 \\ y_4
\end{bmatrix} .
\]
This is equivalent to 
\[
A\begin{bmatrix}
y_3 \\ y_4
\end{bmatrix}
=\begin{bmatrix}
y_1 \\ y_2
\end{bmatrix}.
\]
This implies that the isotropy group of $V$ is isomorphic 
with ${\rm SL}_2\times \Bbb{C}^3$ as a variety (here the 
$\Bbb{C}^3$ is given by the free variables $y_3$, $y_4$ and $g_{35}$) 
which is connected. 
\end{ex}
%

%
\section*{Appendix}
%

Here we include the list of all parabolic subgroups 
in the exceptional Lie algebras where there is a 
Richardson element in the first graded part $\mathfrak{g}_1$ 
and where $G_x=P_x$ holds for Richardson elements.
We list them as in our article~\cite{bw}. In the case 
of ${\rm E}_7$ we omit the three parabolic subgroups 
with $P_x\subsetneq G_x$ (examples $5$, $20$, $25$).

\[
\begin{array}{cc} 
 {\rm G}_2 & {\rm F}_4 \\ 
 & \\ 
(1,1) & (1,1,1,1)\\ 
(1,0) & (1,1,0,1)\\ 
(0,0) & (1,1,0,0) \\ 
& (1,0,0,1) \\ 
& (0,1,0,1) \\ 
& (0,1,0,0) \\ 
& (0,0,0,1) \\ 
& (0,0,0,0) 
\end{array} 
\]

\[ 
\begin{array}{rccc} 
 & {\rm E}_6 & {\rm E}_7 & {\rm E}_8 \\ 
 & & & \\ 
1\ &(1,1,1,1,1,1) & (1,1,1,1,1,1,1) & (1,1,1,1,1,1,1,1) \\ 
2\ &(1, 1, 1, 0, 1, 1) & (1, 1, 1, 0, 1, 1, 1) 
& (1, 1, 1, 0, 1, 1, 1, 1) \\ 
3\ &(1, 1, 1, 0, 1, 0) & (1, 1, 1, 0, 1, 0, 1) 
& (1, 1, 1, 0, 1, 0, 1, 1) \\ 
4\ &(1, 1, 0, 1, 0, 1) & (1, 1, 0, 0, 1, 0, 1) 
& (1, 0, 0, 1, 0, 1, 1, 1) \\ 
5\ &(1, 1, 0, 0, 1, 0) &  
& (1, 0, 0, 1, 0, 1, 0, 1)  \\ 
6\ &(1, 1, 0, 0, 0, 1) & (1, 0, 1, 1, 0, 1, 0) 
&  (1, 0, 0, 1, 0, 0, 1, 1)\\ 
7\ &(1, 1, 0, 0, 0, 0) & (1, 0, 1, 0, 0, 1, 0) 
& (1, 0, 0, 1, 0, 0, 1, 0) \\ 
8\ &(1, 0, 1, 1, 0, 1) & (1, 0, 1, 0, 0, 0, 0) 
& (1, 0, 0, 0, 1, 0, 0, 1) \\ 
9\ &(1, 0, 1, 0, 0, 1) & (1, 0, 0, 1, 0, 1, 1) 
& (1, 0, 0, 0, 0, 1, 1, 1) \\ 
10\ &(1, 0, 1, 0, 0, 0) & (1, 0, 0, 1, 0, 1, 0) 
& (1, 0, 0, 0, 0, 1, 0, 1) \\ 
11\ &(1, 0, 0, 1, 1, 1) & (1, 0, 0, 1, 0, 0, 1) 
& (1, 0, 0, 0, 0, 1, 0, 0) \\ 
12\ &(1, 0, 0, 1, 0, 1) & (1, 0, 0, 0, 1, 0, 0) 
& (1, 0, 0, 0, 0, 0, 1, 1) \\ 
13\ &(1, 0, 0, 1, 0, 0) & (1, 0, 0, 0, 0, 1, 1 ) 
& (1, 0, 0, 0, 0, 0, 1, 0) \\ 
14\ &(1, 0, 0, 0, 1, 1) & (1, 0, 0, 0, 0, 1, 0) 
& (1, 0, 0, 0, 0, 0, 0, 1) \\ 
15\ &(1, 0, 0, 0, 1, 0) & (1, 0, 0, 0, 0, 0, 1) 
& (1, 0, 0, 0, 0, 0, 0, 0) \\ 
16\ &(1, 0, 0, 0, 0, 1) & (1, 0, 0, 0, 0, 0, 0) 
& (0, 1, 0, 0, 0, 0, 0, 1)  \\ 
17\ &(1, 0, 0, 0, 0, 0) & (0, 1, 1, 0, 0, 1, 1) 
& (0, 1, 0, 0, 0, 0, 0, 0) \\ 
18\ &(0, 1, 1, 0, 1, 1) & (0, 1, 0, 0, 0, 0, 0) 
& (0, 0, 1, 0, 0, 0, 1, 0)\\ 
19\ &(0, 1, 1, 0, 0, 1) & (0, 0, 1, 0, 0, 1, 0) 
&   (0, 0, 0, 1, 0, 0, 1, 1) \\ 
20\ &(0, 1, 0, 1, 0, 0) &  
& (0, 0, 0, 1, 0, 0, 1, 0)\\ 
21\ &(0, 1, 0, 0, 0, 1) & (0, 0, 1, 0, 0, 0, 0) 
& (0, 0, 0, 1, 0, 0, 0, 1) \\ 
22\ &(0, 1, 0, 0, 0, 0) & (0, 0, 0, 1, 0, 1, 0) 
& (0, 0, 0, 0, 1, 0, 0, 1) \\ 
23\ &(0, 0, 1, 0, 0, 1) & (0, 0, 0, 1, 0, 0, 1) 
&  (0, 0, 0, 0, 1, 0, 0, 0) \\ 
24\ &(0, 0, 1, 0, 0, 0) & (0, 0, 0, 1, 0, 0, 0) 
& (0, 0, 0, 0, 0, 1, 0, 0)\\ 
25\ &(0, 0, 0, 1, 0, 1) &  
& (0, 0, 0, 0, 0, 0, 1, 1) \\ 
26\ &(0, 0, 0, 1, 0, 0) & (0, 0, 0, 0, 1, 0, 0) 
& (0, 0, 0, 0, 0, 0, 1, 0)  \\ 
27\ &(0, 0, 0, 0, 1, 1) & (0, 0, 0, 0, 0, 1, 0) 
& (0, 0, 0, 0, 0, 0, 0, 1)  \\ 
28\ &(0, 0, 0, 0, 1, 0) & (0, 0, 0, 0, 0, 0, 1) 
&  (0, 0, 0, 0, 0, 0, 0, 0)\\ 
29\ &(0, 0, 0, 0, 0, 1) & (0, 0, 0, 0, 0, 0, 0) 
& \\ 
30\ &(0, 0, 0, 0, 0, 0) &  & 
\end{array}\] 

%
%

\bibliographystyle{amsplain}

\end{document}